\documentclass[12pt]{article}
\usepackage{epsfig}
\usepackage{amsmath,amsthm,amsfonts,amssymb,amscd}

\newcommand{\diff}{\operatorname{Diff}}

\newcommand{\pf}{{\flushleft{\bf Proof: }}}

\newcommand{\ra}{{\rightarrow}}

\def\cC{\mathcal{C}}

\def\cE{\mathcal{E}}

\def\cH{\mathcal{H}}
\def\cR{\mathcal{R}}

\def\cO{\mathcal{O}}
\def\cP{\mathcal{P}}

\newtheorem{ttt}{Theorem}
\newtheorem{lem}{Lemma}
\newtheorem{question}{Question}
\newtheorem{defi}{Definition}

\newtheorem{cor}{Corrollary}
\newtheorem{pro}{Proposition}
\begin{document}
\title {Robust transitivity implies almost robust ergodicity }
\author{\sc Ali Tahzibi
 \thanks{This work began at IMPA and finalized at USP-S\~ao Carlos
supported by FAPESP.}\\
 ICMC-USP S\~ao Carlos, Brazil\\
(e-mail: tahzibi@impa.br)} \maketitle
\begin{abstract}
 In this paper we show the relation between robust
transitivity and robust ergodicity for conservative
diffeomorphisms. In dimension 2 robustly transitive systems are
robustly ergodic. For the three dimensional case, we define {\it
almost robust ergodicity} and prove that generically robustly
transitive systems are almost robustly ergodic, if the Lyapunov
exponents are nonzero. We also show in higher dimensions, that
under some conditions robust transitivity implies robust
ergodicity.
\end{abstract}

\section*{Introduction}

We shall address here the question of how the important concepts
of topological transitivity and ergodicity (metric transitivity)
are related. It is easy to verify that if $ m :=$ Lebesgue measure
is ergodic or $f$ preserves any ergodic probability measure
 which gives a positive mass to
open balls, then $f$ is topologically transitive. That is, there
exists a point $ x \in M $ such that its orbit is dense, or for
any two open sets $ U$ and $V$, there exists $ n \in \mathbb{N} $
such that $ f^n (U ) \cap V \not= \emptyset$.

On the other hand, we know that the converse implication is not
true i.e transitivity is not enough to garante ergodicity. In
fact, by an example of Furstenberg, we even know minimal and
non-ergodic diffeomorphism.

Although transitivity and ergodicity are different notions
(topological and metric), we want to relate them when they persist
in a neighbourhood of a diffeomorphism. To
 obtain ergodic results we need
  more regularity, for instance H\"older continuity
   for the derivative of the diffeomorphism is generally a necessary condition for
    proving ergodicity. It is relevant to remember that even for
    $C^1$ Anosov diffeomorphisms
    preserving volume, the ergodicity is not verified.

    Let us define
     $\diff_m^{1 + } (M ):= \bigcup_{\alpha > 0 } \diff_m^{ 1 + \alpha}$ and prove
     ergodicity results in this set.
\begin{defi}
We say that $f \in \diff_{m}^1 (M)$ is $C^1$-robustly transitive
     (resp.
$C^1$-robustly ergodic), if there exists an open set in $C^1$
topology $ U \subset \diff^1 (M )\\
 ($resp. $ U \subset
\diff^1_{ m } (M )$) such that any $g$ in $U$  (resp. in $U \cap
\diff^{1+} (M)$)is also topologically transitive (resp. ergodic).
\end{defi}
A $Df$-invariant splitting $ E \oplus F $ of
$TM$ is called dominated splitting if the fibers of the bundles
have constant dimension on whole manifold and there is $ \lambda <
1 $ such that:

$$ \| Df | E_x \| . \| Df^{-1} | F_{f (x )} \| \leq \lambda  \quad \text{for all } \quad
x \in M  $$
In~\cite{BDP} the authors show that every $C^1-$robustly
transitive diffeomorphism has dominated splitting. More precisely
for $C^1$-robustly transitive diffeomorphisms,
 $TM = E_1 \oplus E_2 \cdots \oplus E_k$ and this decomposition
is dominated, moreover $Df$ behaves ``hyperbolic" for volume in
$E_1 , E_k$, more precisely for some $C > 0$ :
                    $$ | \det (Df^{ - n } | E^{k} (x ) | \leq C \lambda^n  $$
                    $$ | \det (Df^{  n } | E^{1} (x ) | \leq C \lambda^n   $$

As for surface diffeomorphisms, $E^1$ and $E^k$ are one
dimensional, $C^1$-robust transitivity implies
 a global hyperbolic structure
 or in other word, robustly transitive diffeomorphisms are Anosov.
Robustly transitive diffeomorphisms defined on three-manifolds may
be non-Anosov. In this case the tangent bundle of the ambient
manifold can be split in the following ways:
\begin{itemize}
\item $ TM = E^{cs} \oplus E^u$  where $E^u$ is uniformly expanding one dimensional
and $E^{cs}$ can not be split in whole manifold and it is ``volume
contracting":
$$  | \det (Df^n | E^{cs} (x ) | \leq C \lambda^n  \quad \text{ for } \quad \lambda < 1$$
\item $ TM = E^s \oplus E^{cu}$
\item $ TM = E^s \oplus E^c \oplus E^u $, three subbundles are nontrivial
 and $E^c$ is not uniformly hyperbolic. This case is called strongly partially hyperbolic.

\end{itemize}
Let $\cP\cH^r (M )$ (resp. $\cP\cH^r_m (M )$) be the set of
$C^r$-partially hyperbolic (resp. conservative partially
hyperbolic) diffeomorphisms.
The following questions are of interest :
\begin{question}
Let $f \in \diff^1_{m} (M )$ be $C^1$-robustly transitive. May $f$
be approximated by robustly ergodic diffeomorphisms, or even is it
true that any robustly transitive diffeomorphism is robustly
ergodic ?
\end{question}
In the surface diffeomorphisms case, robust transitivity
diffeomorphisms are Anosov and
 it is well known that  conservative Anosov diffeomorphisms are robustly
 ergodic~\cite{An67}.
 So, there is a positive answer to
 the above question in dimension 2. \\
 For strongly partially hyperbolic diffeomorphisms of three dimensional manifolds,
 Dolgopyat~\cite{Do2} has shown that stably ergodic systems are dense.
 An important thing is to remove zero Lyapunov exponent in the
central direction.

 Although all known examples of robustly
transitive diffeomorphisms have nonzero Lyapunov exponents, it is
more reasonable to ask
 the following:
\begin{question}
May one approximate any robustly transitive diffeomorphism by
another whose Lyapunov exponents are nonzero in a full Lebesgue
measure set.
\end{question}
Absolute continuity of stable and unstable foliations for the
non-uniformly hyperbolic $ C^{ 1 + \alpha}$diffeomorphisms is the
main ingredient to prove
 ergodicity in majority of cases. We remark that a $C^{1 + \alpha}$
 assumption on $f$ is needed to get the absolute continuity of the
 invariant foliations, for an example of non-absolute continuous foliation
  in $C^1$-Anosov case see ~\cite{RY90}.

\begin{ttt} \label{partial}
If $f$ is $C^1$-robustly transitive and $TM = E^u \oplus E^{cs}$
with negative Lyapunov exponents in the $E^{cs}$ direction on a
full Lebesgue measure set (also $C^1$-robustly) then $f$ is
$C^1$-robustly ergodic.
\end{ttt}
The hypothesis about Lyapunov exponents in the above theorem is
also called ``mostly contracting". Mostly expanding is defined
similarly and the theorem is true for them too. The main point in
the above theorem is that for partially hyperbolic diffeomorphisms
the invariant foliations tangent to hyperbolic directions are
continuous.

In Theorem\ref{partial},  all the Lyapunov exponents in the
central direction are of the same sign.  If this does not happen,
that is if there exists directions corresponding to the positive
Lyapunov exponents which do not dominate directions corresponding
to negative Lyapunov exponents, we can not get the same results.
For these cases we prove that any robustly transitive
diffeomorphism of three dimensional manifold is approximated by
{\it almost robustly ergodic} ones. We use a new result of
Bochi-Viana which shows dominated splitting for Oseledets
decomposition of $C^1$ generic conservative diffeomorphisms. Let
us define almost robust ergodicity :
\begin{defi}
$ f \in \diff^{1 + \alpha} (M )$ is said $\epsilon$-ergodic, if
there exists an SRB measure $\rho$ for $f$ with Leb$ (B (\rho ) )
> 1 - \epsilon$ , and in the conservative case $m (C ) > 1 -
\epsilon $, where $C$ is some ergodic component of the Lebesgue
measure.
\end{defi}

\begin{defi}
$f \in \diff^1_m (M ) $ is said  $C^1$-almost robustly ergodic, if
for any $ \epsilon > 0 $ there is an open set $U_{\epsilon}
\subset \diff_m^{1} (M ) $ such that $ g \in U_{\epsilon} \cap
\diff_m^{1 + } $ is $\epsilon$-ergodic.
\end{defi}
By the above definition of almost-robust ergodicity of $f$ we
require that perturbing a little $f$, with a large enough
probability, the Birkhoff averages of continuous function are the
same. In other words the ergodicity near $f$ is `` getting
better". It is clear by definition that if $f \in \diff_m^{1 + }
(M)$ is almost-robustly ergodic, then $f$ is ergodic.
     \begin{ttt} \label{almost}

If $U$ is a $C^1$ open set in $\cP\cH_m^1 (M )$ (dim $M = 3$) of
transitive diffeomorphisms such that for $f \in U$ the central
Lyapunov exponent is nonzero in a full Lebesgue measure set, then
almost-robustly ergodic ones are $C^1$-residual in $U$
     \end{ttt}
As a corollary:
\begin{cor} \label{generic}
 The ergodic diffeomorphisms constitute a generic subset of $U \cap \diff_m^{1 +} $, with the $C^1$ induced topology, if $U \cap \diff_m^{1 +}$ is dense in $U .$
\end{cor}
The density hypothesis in the above corollary is very natural. In
fact, in the symplectic case it is proved that $C^{\infty}$
symplectic diffeomorphisms are
 dense in $C^1$ ones~\cite{Ze77}. For conservative diffeomorphisms it is an open question yet.
 We mention here that another interesting question for robustly
 ergodic diffeomorphisms is the following:
 \begin{question}
 Is it true that any $C^1$-robustly ergodic conservative
 diffeomorphism admits a dominated splitting?
 \end{question}

For more than three dimensional case in~\cite{BoV00}, the authors
give the first example of robustly transitive diffeomorphisms of
$\mathbb{T}^4$ which do not have any hyperbolic direction and
in~\cite{T02} we show the robust ergodicity of such
diffeomorphisms. In fact, the  robust ergodicity of a class of
diffomorphisms of $\mathbb{T}^n$, which do not have any hyperbolic
direction is verified.

\section{Preliminary}

Let $ f: M \rightarrow M $ be a $ C^{ 1 + \alpha}$ diffeomorphism
of a compact  manifold that preserves volume $m$. Oseledets
theorem ~\cite{Os68} states that, for $m-$almost every point $x
\in M $, there exists real numbers $\lambda_1 (x ) > \cdots >
\lambda_{k (x )} (x )  $ and $$ T_x M = E^1_x \oplus \cdots \oplus
E^{ k (x ) }_x$$
such that
           $$ \lim_{ n \rightarrow \infty } \frac{1}{n}\, \log\,
            \| D f^n (x ) (v_j ) \| = \lambda_j (x ) \quad \text{for all }
             \quad v_j \in E^j_x \backslash \{0\} $$
$\lambda_j$'s are Lyapunov exponents which depends measurably on $x$.\\
%
%

%
%
%
Let us mention that for a general $C^{ 1 + \alpha}$ diffeomorphism
preserving an smooth hyperbolic measure (with nonzero Lyapunov
exponents) one has a countable number of ergodic
 components~\cite{Pe77}.  Dolgopyat, Hu and Pesin constructed an example of such diffeomorphism
 with infinitely many open ergodic components. It is obvious that such examples can not be
 topologically
 transitive because open invariant sets are dense in transitive case, so there would exist
 just one ergodic component.\\
We are going to show ergodicity by finding open sets (mod 0) in
the ergodic components of Lebesgue measure. This kind of approach
to prove ergodicity was introduced by Pesin \cite{Pe77}. The idea
is to show local ergodicity (ergodic components are open sets,
a.e) and then by the aid of topological transitivity we get
ergodicity. Observe that to prove ergodicity or equivalently the
existence of a unique ergodic component, we just need to show that
any ergodic component contains an open set (mod 0) and then use
the transitivity as following :

 If $ C_1 , C_2$ are two ergodic component with $U_1 \subset C_1 $, $U_2 \subset C_2 $
 (mod 0), i.e $m ( U_i \backslash C_i) = 0$ then by topological transitivity
  $ \exists n > 0 $ such that
 $ f^n (U_1 ) \cap U_2 \neq \emptyset $. By the invariance of $C_i$
  and the fact that $f$ is
 diffeomorphism and
 and preserve the measure $m$, it comes out that $C_1 \cap C_2 \neq \emptyset$ and $C_1 = C_2.$\\
\begin{proof}(Theorem 1)

 Let $f \in U$ be as in Theorem \ref{partial}, where $U \in \diff_m^1 (M)$
 constitutes of transitive diffeomorphims with non-zero Lyapunov exponents
 in a full Lebesgue measure set.  As $m$ is an invariant measure with non-zero
 Lyapunov exponents,
it has a countable number of ergodic components which are just the
normalization of its restriction to $C_i$  where each $C_i$ is
 a measurable invariant subsets with $m (C_i) > 0$ which we call them
 also by
 ergodic components.

By  continuity of the invariant foliations of partially hyperbolic
diffeomorphisms in \cite{HPS77}, the unstable foliation of $g \in
U $ varies continuously. Using Pesin theory for $g \in \diff^{1
+}_m (M)$, there is a full Lebesgue measure subset $\cR$ such that
 $x \in \cR$ is regular in the following sense.
 $$ \lim_{ n \ra \infty}\frac{1}{n} \sum_{i=0}^{n - 1} \phi
({f^{i}} (x ) ) = \lim_{ n \ra - \infty}\frac{1}{n} \sum_{i=0}^{n
- 1} \phi ({f^{i}} (x ) ) \quad \text{for any $\phi \in C^0 (M)$}
\quad
$$
For any $x \in \cR$ there is a local stable manifold $W^s_{loc}
(x)$ which depend measurably on $x$. Moreover the stable foliation
is absolutely continuous. For any ergodic component $C$, there is
$x \in C$ such that $W_{loc}^s(x ) $ is almost (With respect to
Lebesgue measure of $W^s (x )$) contained in $C \cap \cR $. That
means, there exists $C_x \subset W^s_{\epsilon} \cap C \cap \cR$,
for some $\epsilon > 0$, such that $m_s (W^s_{\epsilon} \backslash
C_x ) = 0$ where $m_s$ is the Lebesgue measure induced on the
stable manifold. Now we saturate $C_x$ by unstable leaves and let
$U_x= \bigcup_{ y \in W_{\epsilon}^s} W^u (y )$. By continuity of
unstable foliation, $U_x$ will contain an open set. On the other
hand by Hopf's argument
 for any $y \in C_x$ and  $z \in W^u (y)$ we have:
 $$ \lim_{ n \ra - \infty}\frac{1}{n} \sum_{i=0}^{n - 1} \phi
({f^{i}} (z ) ) = \lim_{ n \ra - \infty}\frac{1}{n} \sum_{i=0}^{n
- 1} \phi ({f^{i}} (y ) )=$$
$$ \lim_{ n \ra \infty}\frac{1}{n}
\sum_{i=0}^{n - 1} \phi ({f^{i}} (y ) ) = \lim_{ n \ra
\infty}\frac{1}{n} \sum_{i=0}^{n - 1} \phi ({f^{i}} (x ) )$$
Now
by absolute continuity of $W^u$ (unstable foliation) we conclude
that $ \bigcup_{ y \in C_x}
 W^u (y )$ has full measure in $U_x$ and finally as regular points have full measure,
 the ergodic component $C$ contains a total Lebesgue measure subset of the open set $U_x$.
 This is what we want because
 ergodic components are invariant and $f$ is transitive and this
 proves the uniqueness of ergodic components with positive Lebesgue measure
 and consequently ergodicity of the Lebesgue measure.
 \end{proof}

  We mention that, only absolute continuity of $W^u$
  is not enough to get an open set (mod 0) in the above argument,
   but if we have some control from below on
  the size of local unstable manifolds of $y \in W_{local}^s (x )$ then the same argument
  works.
The problem in Theorem \ref{almost} appears when in some ergodic
component with positive measure, a direction corresponding to a
positive  Lyapunov exponent does not dominate the direction
corresponding to the negative one. More precisely let $TM = E^u
\oplus E^{cs}$ and $ \lambda_1 (x ) > \lambda_2 (x )
> \lambda_3 (x ) $ represent the Lyapunov exponents. It may happen
that for an ergodic component $C$ with $m (C) > 0$, $\lambda_2 >
0$ and $E^{cs}$ can not be split into one dimensional dominating
subbundles over $C$.

    We show that the measure of such bad behaved set is small enough for diffeomorphisms
   near to $C^1$-generic diffeomorphisms for which almost everywhere the Oseledets splitting is
   dominated.\\
\section{Dominated splitting of Oseledet splitting}
Now let recall the new result of Bochi-Viana \cite{BoV02} that
gives a $C^1$ generic subset
 of diffeomorphisms for which the Oseledet decomposition is dominated or trivial.
\begin{ttt}[\cite{BoV02}] \label{dichotomy}
For any compact manifold $M$, there exists a $C^1$-residual subset
$\cR$ of $ \diff_m^1 (M )$ such that for every $f \in \cR$, the
Oseledet splitting is dominated or else trivial, at almost every
point.
\end{ttt}
  Let us define $\Lambda_i (f ) = \lambda_1 (f ) + \cdots +
\lambda_i (f )$, where $\lambda_j (f ) = \int_M \lambda_j (x ) dm
$. By semi-continuity arguments it is easy to show  that the
continuity points of $ \Lambda : f \ra (\Lambda_1 (f ) , \cdots ,
\Lambda_{d - 1}
 (f ) ) $ contains a $C^1$-residual subset $ \cR \subset \diff_m^1 (M )$. In \cite{BoV02} the
 authors show that $\forall f \in \cR$,
 the Oseledet splitting is dominated or trivial. We show that diffeomorphisms
 near to this generic subset have greedy ergodic components.
%
%
%
%
%
%
\pf(of Theorem \ref{almost})\\
Let $TM = E^u \oplus E^{cs}$ be the decomposition of tangent
bundle for the diffeomorphisms in $U$.(the other cases are
similar.) Take any $g \in \cR \cap U$, where $\cR$ is the residual
subset in Theorem \ref{dichotomy}. As $g$ is a continuity point of
$\Lambda$ above, for small $\delta_0
> 0$ there exists $\cO$, a $C^1$ neighbourhood of $g$ in
$\diff_m^1 (M )$ such that $ \forall f , f_1 \in \cO$,
 $| \Lambda_i (f_1 ) - \Lambda_i (f ) | < \delta_0 $.\\
As Lyapunov exponents of $f \in U$ are nonzero then by
\cite{Pe77}, $m$ has a countable number
 of ergodic components $C_n$ with $m (C_n ) > 0 $ and $m (\bigcup C_n ) = 1$.\\
Now let $f \in \cO \cap \diff_m^{1 +} (M )$; we will show that $f$
has a large ergodic component in the measure theoretic sense.

 Any ergodic component $C$ is an invariant set and
$m_C$:= the normalized restriction of the Lebesgue measure on $C$
is an ergodic measure for the restriction of $f$ on $C$.
Neglecting zero measure subsets which are irrelevant for our
purposes, we may substitute $C$ by $Supp (m_C)$. Now we claim that
there exists $\bar{C} \subseteq C$ such that
 $m_C ( C \backslash \bar{C})= 0$ and any $x \in \bar{C}$ has a dense orbit in $C$.
 To show the above claim remember that by the ergodicity of
 $m_C$, the basin of $m_C$ has full measure. The basin
 consists of points $z$ such that $\frac{1}{n} \sum_{i=0}^{n-1}
 \delta_{z}$ converges to $m_C$. As $m_C$ has a full support,
   the orbit of the points in the
 basin will be dense in $C$.  Now it is enough to take $\bar C =$
 Basin of
  $m_C$ and the claim is proved. In the following we define two
  kinds of ergodic components:
\begin{enumerate}
\item  $\cC_{good} =$ Union of ergodic components $C_i$ such that
 one of the followings happens:
 \begin{itemize}
 \item $\lambda_2 (x ) < 0 $ for $x \in C_i$
 \item $E_2 \oplus E_3$ is a dominated
splitting on $C_i$ and $\lambda_2 (x ) > 0.$
\end{itemize}
of course for each ergodic component at most one of the
 above items can be satisfied.
 \item $\cC_{bad}$= Complement of $\cC_{good}$
\end{enumerate}

By the arguments in the previous section in the proof of Theorem
\ref{partial} we deduce that all $ C_i \subset \cC_{good} $
contains open set (mod 0), because:
\begin{itemize}
\item If $\lambda_2 (x ) < 0 $ then we are in the setting of
Theorem \ref{partial}. Firstly we use absolute continuity of local
stable manifolds corresponding to $\lambda_2$ and $\lambda_3$ and
then continuity of unstable foliation tangent to $E^u$.
\item If $\lambda_2 (x)= \lambda > 0$ and $E_2$ dominates $E_3$,
then on the corresponding ergodic component $C_i$, the tangent
bundle has the following invariant dominated decomposition: $T_C
M= E^{cu} \oplus E^{cs}$ with $\dim (E^{cs}) = 1$ and $E^{cu} =
E^u \oplus E^2$ where $E^2$ corresponds to $\lambda_2$.
 Now as $f$ is
conservative it is easy to show that the above dominated splitting
has the volume hyperbolicity property: $|\det ( Df^n| E^{cs} (x)|
\leq C \lambda^n $ for some $C > 0 , \lambda < 1$ for any $x \in
C_i$ (see \cite{BDP} for volume hyperbolicity  of dominated
splitting in the coservative case). As $E^{cs}$ is one
dimensional, the volume hyperbolicity is equivalent to uniform
hyperbolicity. Now apply the same argument of the Theorem
\ref{partial} using the absolute continuity of the local unstable
manifolds for the points in $C_i$ and the existence of large
stable manifolds which is a result of uniform hyperbolicity of
$E^{cs}$.
\end{itemize}
 Consequently by  topological transitivity, there is just one egodic component in
 $\cC_{good}$.
For the bad components we prove that they occupy an small region
 in the measure theoretic sense.
  Let us fix some notations. By an $m-$dominated splitting of
 $E_2 \oplus E_3$
 along the orbit of a point $x$ we require that for all $n \in \mathbb{Z}$:
 $$\displaystyle \frac{\| Df_{f^n (x)}^m\mid E_3 \|}{{\bf m} (Df_{f^n(x)}^m \mid E_2)}
  \leq \frac{1}{2}$$
  where ${\bf m}(A)= \| A^{-1}\|^{-1}$.
  By an $m-$dominated splitting over an invariant set $C$ we mean
  $m-$dominated splitting for all orbits in $C$. Observe that
  $m-$dominated splitting extends to the closure of an invariant
  subset.

 Let $\Gamma (f, m)$ denotes the subset of points such that $E_2 \oplus E_3$
does not admit an $m-$dominated splitting and let $\Gamma (f ,
\infty)= \cap_{m \in \mathbb{N}} \Gamma (f, m)$. For the ergodic
components  $ C \subseteq \cC_{bad}$, $\lambda_2 > 0$ and
 $E_2 \oplus E_3$ does not admit an $m-$dominated splitting over
 $C$ for any $m \in \mathbb{N}$. We claim that, there exits a full $m_C$
 measure subset of $C$ which we denote it also by $C$, such that
  $C \subseteq \Gamma(f, \infty)$. To prove the claim
 observe that as $m-$dominated splitting passes to closure of sets,
  for any point $ x \in \bar{C}$ ($z \in \bar{C}$ has dense orbit)
 there is not any $m-$dominated
 splitting of $E_2 \oplus E_3$ along the orbit of $x$. This shows
 that
 $\cC_{bad}\subseteq \Gamma(f, \infty)$ (mod 0).

Denote $$J (f) = \displaystyle \int_{\Gamma(f, \infty)}
\displaystyle \frac{\lambda_2 - \lambda_3}{2} dm (x)$$
 We apply
the following proposition of \cite[Proposition 4.17]{BoV02}:
 \begin{pro}
 Given any $\delta > 0$ and $\epsilon > 0$, there exists a
 diffeomorphism $f_1$, $\epsilon$ near to $f$ such that
 $$
 \int_M \Lambda_2 (f_1 , x ) dm (x) < \int_M \Lambda_2 (f , x ) dm (x)-
 J (f) + \delta $$
 \end{pro}
By the above proposition we will deduce that if the measure of bad
components is not small enough then after perturbing $f$ a little,
the average of $\lambda_1 + \lambda_2$ drastically drops.

As $\cC_{bad}\subseteq \Gamma(f, \infty)$ and on $ \cC_{bad} ,
\lambda_2 (x ) > 0 $ by the above proposition it turns out that
   $$ \Lambda_2 (f ) - \Lambda_2 (f_1 ) \geq  \frac{1}{2} \int_{\cC_{bad}}
    (\lambda_2 - \lambda_3 ) (f ) dm - \delta \geq  \frac{1}{2} \int_{\cC_{bad}}
     - \lambda_3  (f ) dm - \delta$$
$$ \geq m (\cC_{bad} ) \, inf_{ x \in \cC_{bad}} \frac{- \lambda_3 (f , x )} {2}
- \delta $$
Now we observe that the infimum above is bounded away from zero
uniformly in $\cO$.
\begin{lem}
There is some $\alpha < 1$, depending only to $g$ such that:
$\lambda_3 (x ) \leq \log (\alpha ) $ for all $x \in \cC_{bad}$.
\end{lem}
\begin{proof}
As $g$ is volume hyperbolic and partially hyperbolic $TM = E^{u}
\oplus E^{cs}$, so $ \det (Df| E^{cs} (x ) ) < \alpha < 1 $ for
all $x \in M $ and we can take $\alpha$ uniform in a $C^1$
neighborhood of $g$, just because $g \ra E^{cs} (x , g )$ is
continuous for $C^1$ topology in the partially hyperbolic space.
Take $x \in \cC_{bad}$ then:
$$ \lambda_2 (x ) + \lambda_3 (x ) =
 \lim_{n \ra \infty} \frac{1}{n} \sum_{i = 0}^{n - 1} \log\, \det
  (Df | E^{cs} (f^i (x ) ) ) \leq \log (\alpha )$$ and as
  $\lambda_2 (x ) > 0 $ we have $  \lambda_3 (x ) < \log (\alpha )$.
\end{proof}
So by the above lemma:
$$ inf_{x \in \cC_{bad}} \frac{- \lambda_3 (f , x )}{2} \geq \frac{ - \log (\alpha )}{2} $$
and
$$ \delta_0 \geq \Lambda_2 (f ) - \Lambda_2 (f_1 ) \geq m (\cC_{bad} )
 \frac{ - \log (\alpha )}{2} - \delta$$
 We get then $ m (\cC_{bad} ) \leq \frac{\displaystyle \delta + \delta_0}
 { \displaystyle - \log (\alpha )} $. Taking $ \delta_0$ and $\delta$
  small, for $f \in \cO \cap \diff^{1+}_m (M)
 ,  m (\cC_{bad} )$ is small enough and so $m (\cC_{good})$ is large and ``{\it ergodicity is
 getting better"}.\\

Finally we conclude The Corollary ~\ref{generic} as follows.
Taking $\cE_n:= \frac{1}{n}$-ergodic diffeomorphisms in $
\diff_m^{1 +} \cap U $, then $\cE_n$ is open and dense in the
$C^1$ induced topology , so $\cE = \bigcap \cE_n$ is a residual
subset and $f \in \cE$ is ergodic.



\bibliographystyle{plain}
\bibliography{bib}
\end{document}